# THE FUNDAMENTAL LEMMA FOR
# THE SHALIKA SUBGROUP OF $GL(4)$

BY SOLOMON FRIEDBERG AND HERVÉ JACQUET

University of California, Santa Cruz
Columbia University

DEPARTMENT OF MATHEMATICS, UNIVERSITY OF CALIFORNIA SANTA CRUZ, SANTA CRUZ, CA 95064
*E-mail address*: friedbe@cats.ucsc.edu

DEPARTMENT OF MATHEMATICS, COLUMBIA UNIVERSITY, NEW YORK, NY 10027
*E-mail address*: hj@math.columbia.edu

Typeset by $\mathcal{A}_{\mathcal{M}}\mathcal{S}$-TEX

# ABSTRACT

The authors establish the fundamental lemma for a relative trace formula. This trace formula compares generic automorphic representations of $GSp(4)$ with automorphic representations of $GL(4)$ which are distinguished with respect to a character of the Shalika subgroup, the subgroup of matrices of $2 \times 2$ block form $\begin{pmatrix} g & 0 \\ X & g \end{pmatrix}$. The fundamental lemma, giving the equality of two orbital integrals, amounts to a comparison of certain exponential sums arising from these two different groups. On the $GSp(4)$ side one has the Kloosterman sums for this group, and on the $GL(4)$ side certain, new, relative Kloosterman sums. To show that these are equal for each relevant Weyl group element, the authors compute the Mellin transforms of the sums and match them in all cases.

1991 *Mathematics Subject Classification.* Primary 11F70; Secondary 11F46, 11F72, 11L05, 22E50.

Research supported by National Science Foundation grants DMS-8821762 and DMS-9123845 (Friedberg), and DMS-8801759 and DMS-9101637 (Jacquet). Friedberg's research also supported by the Sloan Foundation. Research at MSRI is supported in part by National Science Foundation grant DMS-9022140.

# CONTENTS

CHAPTER I

INTRODUCTION AND STATEMENT OF RESULTS

1. Introduction

(1.1) Let $K$ be a number field and let $G$ be the group $GL(4)$ regarded as an algebraic group over $K$. Suppose that $\pi$ is a cuspidal automorphic representation of $G$. Thus, conjecturally, $\pi$ is associated to a representation $\tau$ of the hypothetical Langlands group $L_K$ into $^LG^0 = GL(4, \mathbf{C})$. Suppose furthermore that $\pi$ is self-contragredient. Then either the symmetric square $L$-function or the exterior square $L$-function of $\pi$ has a pole at $s = 1$. Assume that the exterior square $L$-function has a pole at $s = 1$. Conjecturally, this means that $\wedge^2 \tau$ has a fixed vector, and hence that $\tau$ factors through the symplectic group in 4 variables $Sp(4, \mathbf{C})$, which is the $L$-group of $SO(5, K)$. Thus it is predicted that $\pi$ is the functorial image of a cuspidal automorphic representation $\pi'$ (or rather a packet of representations) of the group $SO(5)$. On the other hand, the condition that the exterior square $L$-function has a pole at $s = 1$ can be read directly on the space of the representation $\pi$, due to the integral representation for the exterior square $L$-function given by Jacquet and Shalika [JS]. More precisely, let $H$ be the Shalika subgroup of $GL(4)$

$$H = \left\{ \begin{pmatrix} g & \\ & g \end{pmatrix} \begin{pmatrix} I_2 & \\ X & I_2 \end{pmatrix} \,\bigg|\, g \in GL(2), X \in M_{2\times 2} \right\},$$

again regarded as an algebraic group over $K$. If $\psi$ is an additive character of the adeles $\mathbb{A}$ of $K$ trivial on $K$ define a character still denoted $\psi$ of $H(\mathbb{A})$ by

(1.1) $$\begin{pmatrix} g & \\ & g \end{pmatrix} \begin{pmatrix} I_2 & \\ X & I_2 \end{pmatrix} \mapsto \psi(\mathrm{tr} X).$$

Then the exterior square $L$-function has a pole at $s = 1$ if and only if $\pi$ is *distinguished* with respect to $\psi$, in the sense that there is a vector $\phi$ in the representation $\pi$ such that

$$\int \phi(h)\, \psi(h)\, dh \neq 0.$$

More generally, we can consider automorphic representations which are essentially self-contragredient, that is, self-contragredient up to twist by a character. This amounts to replacing the exterior square $L$-function in the above discussion by a suitable twist. A pole at $s = 1$ conjecturally implies that the representation $\tau$ factors through the group of symplectic similitudes $GSp(4, \mathbf{C})$, which is the $L$-group of $GSp(4, K)$. This means that the representation $\pi$ should be distinguished with respect to a character of the form $h \mapsto \chi(\det g)\, \psi(\mathrm{tr} X)$.

(1.2) The above considerations suggest the existence of a trace formula of the following type. For simplicity, we restrict ourselves to the study of self-contragredient representations with trivial central character, and we consider only the character $\psi$ of the Shalika subgroup. Let $\Phi$ be a smooth function of compact support on $G(\mathbb{A})$. Define as usual a kernel function

$$K_\Phi(g_1, g_2) = \sum_{\xi \in G(K)} \Phi(g_1^{-1} \xi g_2).$$

Similarly, let $G'$ be the group $GSp(4)$ regarded as an algebraic group over $K$, let $\Phi'$ be a smooth function of compact support on $G'(\mathbb{A})$, and define the kernel function

$$K_{\Phi'}(g'_1, g'_2) = \sum_{\xi' \in G'(K)} \Phi'(g_1'^{-1} \xi' g'_2).$$

The trace formula in question takes the form:

$$(1.2) \quad \int K_\Phi(h^{-1}, n)\, \psi(h)\, \theta(n)\, dh\, dn = \int K_{\Phi'}({}^t n_1, n_2 z)\, \theta'(n_1^{-1} n_2)\, dn_1\, dn_2\, dz.$$

Here $h$ is integrated over the Shalika subgroup, $n, n_1, n_2$ are integrated over the unipotent radicals of the respective Borel subgroups, $\theta, \theta'$ are nondegenerate characters, and $z$ is integrated over the center of $GSp(4)$. Both sides can be computed in terms of orbital integrals, that is, integrals taken over double cosets of the appropriate groups. Not all double cosets contribute. There is a bijection between the relevant double cosets of the two groups, and the above identity holds if the corresponding orbitals integrals are equal. This is described in more detail below. Assuming that there are enough pairs $\Phi, \Phi'$ satisfying this condition of matching orbital integrals, spectral expansions of the kernels should allow one to derive from this identity another identity asserting that

$$\int K^{\text{cusp}}(h^{-1}, n)\, \psi(h)\, \theta(n)\, dh\, dn = \int K^{\text{cusp}}({}^t n_1, n_2 z)\, \theta'(n_1^{-1} n_2)\, dn_1\, dn_2\, dz.$$

Here the exponent cusp is used to indicate that only the cuspidal part of the kernel is considered. (In fact, one expects that some discrete terms will have to be added to the cuspidal kernel.)

What are the applications of such an identity? Via the twisted trace formula, one expects to be able to establish that (essentially) self-contragredient automorphic representations of $GL(4)$ whose (possibly twisted) exterior square $L$-function has a pole are functorial images from the group $GSp(4)$. The present trace formula would then establish that in any global packet of representations of $GSp(4)$ there is at least one element with a global Whittaker model.

We caution that the details of this program are far from completed. In particular, we have not discussed the relation between representations of $SO(5)$ and $GL(4)$. In the case at hand, the Weil representation might be a quicker way to arrive at the same result. However, for higher ranks, the Weil representation is not available and one has to resort to a trace formula of the above type. In that case the group $GSp(4)$ is to be replaced by a spin group whose $L$-group is the $GSp(2n, \mathbf{C})$.

(1.3) If $\Phi$ and $\Phi'$ are each products of local functions $\Phi_v$, $\Phi'_v$, the matching of orbital integrals for (1.2) reduces to a local question. At almost all places the functions $\Phi_v$, $\Phi'_v$ are the characteristic functions of the respective maximal compact subgroups. In this paper, we match the orbital integrals of these characteristic functions, that is, we establish the "fundamental lemma" for the trace formula at hand. The orbital integral on $GSp(4)$ is of Kloosterman type, while the orbital integral on $GL(4)$ is a new sort of relative Kloosterman integral. This matching is the first crucial step in establishing the formula (1.2). However, our method of proof is rather computational, and thus leaves something to be desired. Essentially, we have to show that functions of two (or one) variables are equal. We do so by showing that their formal Mellin transforms are equal, and this is a series of lengthy calculations. It is remarkable that the resulting formulae are relatively simple. Thus, it is possible that the Mellin transforms of our functions have some independent interest.

This use of formal Mellin transforms was introduced, in the context of quadratic base change on $GL(2)$, by Jacquet and Ye [JY1]. A number of additional papers related to the relative trace formula are listed in the References.

The main theorem of this work is given below in Section 4 of this Chapter. This theorem was announced in [FJ].

(1.4) In the rest of this paper we pass to the local situation. Throughout we fix the following notation. Let $K$ now denote a nonarchimedean local field of residual characteristic not equal to 2, $\mathcal{O}$ be its ring of integers, and $\varpi$ be a local uniformizer of $K$. Let $q$ denote the cardinality of the residue field $\mathcal{O}/\varpi\mathcal{O}$, and $|\ |$ denote the normalized absolute value on $K$, so that $|\varpi| = q^{-1}$. Fix an additive character $\psi$ of $K$ with the property that the largest fractional ideal on which $\psi$ is trivial is precisely $\mathcal{O}$.

(1.5) Equations are numbered consecutively throughout each section. The fifth numbered equation of Chapter 1, Section 2 is referred to as (2.5) throughout Chapter 1, and as (1.2.5) in the other chapters. The numbering of Lemmas, Propositions, and Theorems (together) is similar. Sections are also divided into numbered paragraphs, for the convenience of the reader.

## 2. The Shalika Subgroup $H$ of $GL(4)$: Relevant Double Cosets and Kloosterman Integrals

(2.1) For convenience, let us write $H$ instead of $H(K)$ for the Shalika subgroup of $GL(4, K)$:

$$H = \left\{ \begin{pmatrix} g & \\ & g \end{pmatrix} \begin{pmatrix} I_2 & \\ X & I_2 \end{pmatrix} \,\bigg|\, g \in GL(2,K), X \in M_{2\times 2}(K) \right\}.$$

If $B$ is a Borel subgroup of $GL(2, K)$, let $H(B)$ be the subgroup of the Shalika subgroup defined by

$$H(B) = \left\{ \begin{pmatrix} b & \\ & b \end{pmatrix} \begin{pmatrix} I_2 & \\ X & I_2 \end{pmatrix} \,\bigg|\, b \in B, X \in M_{2\times 2}(K) \right\}.$$

When working with the Shalika subgroup (this Section and Chapters 2 and 3), we let $N$ denote the subgroup of $GL(4,K)$ consisting of upper triangular unipotent matrices.

The additive character $\psi$ gives rise to a character of the group $H$, which we shall again denote by $\psi$, by the formula (1.1). In addition, let $\theta$ denote the character of $N$ given by
$$\theta(n) = \psi(n_{1,2} + n_{2,3} + n_{3,4}).$$

(2.2) A double coset $HrN$ of $GL(4,K)$ is said to be $H$-relevant if the map
$$hrn \mapsto \psi(h)\theta(n)$$
is well defined. This is the case if and only if

(2.1) $$\theta(n) = \psi(rnr^{-1}) \qquad \text{for all } n \in r^{-1}Hr \cap N.$$

It is precisely these double cosets which enter into the relative trace formula.

To make these explicit, let $W$ denote the Weyl group of $GL(4)$ with respect to the standard Borel subgroup $B_4$. Then $W$ is isomorphic to the permutation group $S_4$. To make this explicit, associate with each permutation $\sigma \in S_4$ the permutation matrix with entries 1 in the $(i, \sigma(i))$ positions and 0 elsewhere. Note that with this convention, permutations are multiplied from the *right*, so $(12)(13)=(123)$ rather than $(132)$, and if $w$ is the permutation matrix corresponding to $\sigma$, then $w^{-1}m_{i,j}w = m_{\sigma(i),\sigma(j)}$. Let diag( ) denote the $4 \times 4$ diagonal matrix. Then we have

PROPOSITION 2.1. *The $H$-relevant double cosets of $GL(4,K)$ are*
1. $H\,(1243)\,N$
2. $H\,\text{diag}(-a_1, a_1, 1, 1)\,N$
3. $H\,(132)\,\text{diag}(a_1, 1, 1, 1)\,N$
4. $H\,(12)\,\text{diag}(a_1 a_2, a_1, 1, 1)\,N$

*with $a_1, a_2 \in K^\times$.*

PROOF. We first determine a complete set of coset representatives for the double cosets $H \backslash GL(4,K) / B_4$. Let $P_{2,2}$ be the standard parabolic subgroup of type (2,2), and let $W_{2,2}$ be the subgroup of $W$ corresponding to $P_{2,2}$. Identify $W$ with the group of $4 \times 4$ permutation matrices as above. By the Bruhat decomposition, we have the disjoint union
$$GL(4,K) = \bigcup_{w \in W_{2,2} \backslash W} P_{2,2}\, w\, B_4.$$

Now a set of coset representatives for $H \backslash P_{2,2}$ is given by matrices of the form $\begin{pmatrix} g & \\ & I_2 \end{pmatrix}$ with $g \in GL(2,K)$. Such a $g$ may be written $g = b_1 w_1 b_2$ with $w_1$ a $2 \times 2$ permutation matrix and $b_1$ and $b_2$ either upper or lower triangular (with independent choices). Changing on the left by an element of the Shalika subgroup,
$$H \begin{pmatrix} g & \\ & I_2 \end{pmatrix} w\, B_4 = H \begin{pmatrix} w_1 b_2 & \\ & b_1^{-1} \end{pmatrix} w\, B_4.$$

But for correct choices of upper/lower triangularity of the $b_i$, they may each be pulled past the $w$ and absorbed in $B_4$. It follows that for any $w \in W$,

$$P_{2,2} w B_4 = H w B_4 \cup H(12) w B_4.$$

This gives the decomposition

(2.2) $$GL(4, K) = \bigcup_{<(34)>\backslash W} H w B_4. \quad \text{(disjoint union)}$$

Next let us determine which of the twelve double cosets $HwB_4$ in (2.2) may contain a relevant double coset $HwaN$, for some diagonal $a$. Suppose that $w$ corresponds to a permutation $\sigma$. Consideration of the characters $\psi$ and $\theta$ shows that (2.1) holds only when the following conditions are satisfied:

(1) If $\sigma(3) < \sigma(1)$, then $\sigma(1) = \sigma(3) + 1$.
(2) If $\sigma(4) < \sigma(2)$, then $\sigma(2) = \sigma(4) + 1$.
(3) If $\sigma(3) < \sigma(2)$, then $\sigma(2) \neq \sigma(3) + 1$.
(4) If $\sigma(4) < \sigma(1)$, then $\sigma(1) \neq \sigma(4) + 1$.
(5) If $\sigma(1) < \sigma(2)$ and $\sigma(3) < \sigma(4)$, then either $\sigma(1) + 1 < \sigma(2)$ and $\sigma(3) + 1 < \sigma(4)$, or $\sigma(1) + 1 = \sigma(2)$ and $\sigma(3) + 1 = \sigma(4)$.

Using these conditions, one finds that if a double coset $HwaN$ is relevant, then (modulo left multiplication by (34)), $\sigma$ must be one of the permutations 1, (12), (132), and (1243).

To complete the proof of the Proposition, one checks that for each of these four permutations, the double cosets of the form $HwaN$ satisfying (2.1) are precisely those of the form listed. We omit this calculation.

(2.3) For each relevant double coset, we now introduce an integral of the unit element $\Xi$ of the $GL(4)$ Hecke algebra. The integral is to be compared with an appropriate symplectic integral. To facilitate this comparison, it is useful to consider a Mellin transform of the integral against a quasicharacter of the appropriate torus.

The integral is described as follows. Let $\Xi$ be the characteristic function of the maximal compact subgroup $\mathcal{K} = GL(4, \mathcal{O})$ of $GL(4, K)$. Suppose that the double coset $HwaN$ is relevant. Then we consider

(2.3) $$F(w, a) = \int_N \int_{H/H \cap wNw^{-1}} \Xi(hwan) \, \psi(h) \, \theta(n) \, dh \, dn,$$

where $dh$ is the quotient measure obtained from the normalized right Haar measure on $H$.

A first simplification of (2.3) is given as follows. Let $B$ denote a Borel subgroup of $GL(2, K)$, and let $\mathcal{K}_1$ denote the compact subgroup $GL(2, \mathcal{O})$. Then with suitable normalizations of the Haar measures

$$\int_{GL(2,K)} dg = \int_{\mathcal{K}_1} dk_1 \int_B db.$$

Moving the matrix $\begin{pmatrix} k_1 & 0 \\ 0 & k_1 \end{pmatrix} \in \mathcal{K}$ to the left and using the $\mathcal{K}$-invariance of $\Xi$, we may do the $\mathcal{K}_1$ integral. We find that the integral (2.3) is unchanged if we replace $H$ by the subgroup $H(B)$.

We shall compute the integral $F((1243), I_4)$ directly momentarily. To study the remaining integrals, let $\chi_1, \chi_2$ be quasicharacters of $K$. Let us introduce the following formal Mellin transforms of $F(w, a)$:

$$(2.4) \qquad M_F(\chi_1; 1) = \int_{K^\times} F(1, \operatorname{diag}(-a_1, a_1, 1, 1))\, \chi_1(a_1)\, d^\times a_1$$

$$(2.5) \qquad M_F(\chi_1; (132)) = \int_{K^\times} F((132), \operatorname{diag}(a_1, 1, 1, 1))\, \chi_1(a_1)\, d^\times a_1$$

$$(2.6) \qquad M_F(\chi_1, \chi_2; (12)) = \int_{(K^\times)^2} F((12), \operatorname{diag}(a_1 a_2, a_1, 1, 1))\, \chi_1(a_1)\, \chi_2(a_2)\, d^\times a_1\, d^\times a_2.$$

The meaning of these formal integrals is the following. Let $\phi_A$ be the characteristic function of the compact set $q^{-A} \leq |x| \leq q^A$, and introduce the integrals

$$(2.7) \qquad M_F(\chi_1; 1; A) = \int_{K^\times} F(1, \operatorname{diag}(-a_1, a_1, 1, 1))\, \chi_1(a_1)\, \phi_A(a_1)\, d^\times a_1,$$

$$(2.8) \qquad M_F(\chi_1; (132); A) = \int_{K^\times} F((132), \operatorname{diag}(a_1, 1, 1, 1))\, \chi_1(a_1)\, \phi_A(a_1)\, d^\times a_1,$$

$$M_F(\chi_1, \chi_2; (12); A) = \int_{(K^\times)^2} F((12), \operatorname{diag}(a_1 a_2, a_1, 1, 1))\, \chi_1(a_1)\, \chi_2(a_2)$$
$$(2.9) \qquad\qquad\qquad\qquad \times \phi_A(a_1)\, \phi_A(a_2)\, d^\times a_1\, d^\times a_2.$$

The integrands then have compact support. If the quasicharacters $\chi_1$, $\chi_2$ are ramified, then these integrals have limits as $A$ tends to infinity. The formal Mellin transform is by definition this limit. In (2.7) and (2.8), when the quasicharacter $\chi_1$ is unramified, these integrals represent Laurent polynomials in $\chi_1(\varpi)$. As $A$ tends to infinity, these Laurent polynomials approach a Laurent series; this is by definition the formal Mellin transform. Here the topology on the space of Laurent series is that of pointwise convergence of the coefficients. The treatment of (2.9) is similar; the formal Mellin transform is a Laurent series in $\chi_1(\varpi)$ and $\chi_2(\varpi)$ when both quasicharacters are unramified, and is a Laurent series in one variable in the mixed cases.

For a fixed relevant $w$, we occasionally suppress the dependence of $M_F$ on $w$.

(2.4) The quantity $M_F(\chi_1, \chi_2; (12))$ is evaluated in Chapter 2 of these notes. In Chapter 3 we evaluate the remaining two Mellin transforms. First, though, let us dispose of $F((1243), I_4)$.

PROPOSITION 2.2. $F((1243), I_4) = 1$.

PROOF. A direct computation shows that with $w = (1243)$

$$H \cap wNw^{-1} = \left\{ \begin{pmatrix} 1 & x & & \\ & 1 & & \\ & & 1 & x \\ & & & 1 \end{pmatrix} \begin{pmatrix} 1 & & & \\ & 1 & & \\ r & s & 1 & \\ & u & & 1 \end{pmatrix} \right\}.$$

Replacing $H$ by $H(B)$ where $B$ is the standard Borel subgroup of $GL(2)$ consisting of upper triangular matrices, we thus arrive at the expression

$$F((1243), I_4) = \int \Xi \left( \begin{pmatrix} 1 & & & \\ & 1 & & \\ & & 1 & \\ \gamma & & & 1 \end{pmatrix} \begin{pmatrix} t_3 & & & \\ & t_4 & & \\ & & t_3 & \\ & & & t_4 \end{pmatrix} (1243)n \right)$$

$$\times \theta(n) \left| \frac{t_3}{t_4} \right|^2 d\gamma\, d^\times t_3\, d^\times t_4\, dn$$

$$= \int \Xi \left( \begin{pmatrix} t_3 & & & \\ & t_3 & & \\ & \gamma t_3 & t_4 & \\ & & & t_4 \end{pmatrix} n \right) \theta(n) \left| \frac{t_3}{t_4} \right|^2 d\gamma\, d^\times t_3\, d^\times t_4\, dn,$$

where for the last step we have moved (1243) to the left and used the $\mathcal{K}$-invariance of $\Xi$. But a matrix calculation (or the use of Plücker coordinates, as in Chapter 2 below) shows that the matrix

$$\begin{pmatrix} t_3 & & & \\ & t_3 & & \\ & \gamma t_3 & t_4 & \\ & & & t_4 \end{pmatrix}$$

is in $\mathcal{K}N$ if and only if it is in $\mathcal{K}$. Thus the integrand is nonzero only when $n \in \mathcal{K}$, and it is easy to see that the integral reduces to 1.

## 3. The Symplectic Group $GSp(4)$: Relevant Double Cosets and Kloosterman Integrals

(3.1) We consider the similitude symplectic group $G'(K) = GSp(4, K)$ for the skew symmetric matrix

(3.1) $$w_0 = \begin{pmatrix} 0 & w_s \\ -w_s & 0 \end{pmatrix}$$

where we have set:

$$w_s = \begin{pmatrix} 0 & 1 \\ 1 & 0 \end{pmatrix}.$$

We denote by $T$ the group of diagonal matrices in $G'(K)$. It consists of all matrices of the form

(3.2) $$t = \mathrm{diag}(a_1, a_2, a_2^{-1}\lambda, a_1^{-1}\lambda).$$

The scalar $\lambda$ is the similitude ratio. If $\lambda = 1$ the matrix $t$ is in the symplectic group and the positive roots take the following values on $t$:

(3.3) $$\alpha_1 = a_1 a_2^{-1},\ \alpha_2 = a_2^2,\ \alpha_1\alpha_2 = a_1 a_2,\ \alpha_1^2 \alpha_2 = a_1^2.$$

We let $Z$ be the center of $G'(K)$. When working with the group $GSp(4, K)$ (this Section and Chapters 4 and 5), we let $N$ denote the group of upper triangular *symplectic* matrices with unit diagonal. Thus $N$ is the semi-direct product of the group $V$ of matrices of the form:

(3.4) $$v = \begin{pmatrix} 1 & x & 0 & 0 \\ 0 & 1 & 0 & 0 \\ 0 & 0 & 1 & -x \\ 0 & 0 & 0 & 1 \end{pmatrix},$$

and the group $U$ of matrices of the form:

(3.5) $$u = \begin{pmatrix} 1 & 0 & p & r \\ 0 & 1 & s & p \\ 0 & 0 & 1 & 0 \\ 0 & 0 & 0 & 1 \end{pmatrix}.$$

We define a generic character $\theta$ of $N$ by

(3.6) $$\theta(u) = \psi(s),\ \theta(v) = \psi(x).$$

Also, let $\overline{N}$ be the opposite (transposed) subgroup to $N$.

(3.2) A double coset $\overline{N}rN$ is said to be relevant if the map

$$^t n_1 r n_2 \mapsto \theta(n_1 n_2)$$

is well defined. This is the case if and only if

(3.7) $$\theta(n) = \theta(^t(rnr^{-1})) \qquad \text{for all } n \in r^{-1}\overline{N}r \cap N.$$

It is precisely these cosets which enter into the relative trace formula. Relevant double cosets are easily described: as in any split group, they are associated with the standard parabolic subgroups $P = MU$ of $G'$. Let $w_M$ be the longest element in $W \cap M$ and $t$ an element of the center of $M$. Then the double coset of $w_M t$ is relevant and all relevant double cosets have a unique representative of this form for a suitable $M$.

More precisely, let $w_1$ and $w_2$ be the following matrices:

(3.8) $$w_1 = \begin{pmatrix} 1 & 0 & 0 & 0 \\ 0 & 0 & -1 & 0 \\ 0 & 1 & 0 & 0 \\ 0 & 0 & 0 & 1 \end{pmatrix},\quad w_2 = \begin{pmatrix} 0 & 1 & 0 & 0 \\ 1 & 0 & 0 & 0 \\ 0 & 0 & 0 & 1 \\ 0 & 0 & 1 & 0 \end{pmatrix}.$$

Then the a set of representatives for the relevant double cosets consists of the matrices of the form

(1) $t$
(2) $w_1 \operatorname{diag}(a_1, a_2, -a_2, -a_1^{-1}a_2^2)$
(3) $w_2 \operatorname{diag}(a_1, a_1, a_1^{-1}\lambda, a_1^{-1}\lambda)$
(4) $w_0 \operatorname{diag}(a_1, a_1, a_1, a_1)$.

(3.3) Let $\Xi'$ be the characteristic function of $GL(4, \mathcal{O}) \cap GSp(4, F)$. For each relevant double coset $\overline{N}wtN$ we consider the integrals

$$(3.9) \qquad J(w,t) = \int_{N \cap w^{-1}\overline{N}w \backslash N} \int_N \Xi'({}^t n_2 w t n_1) \theta(n_1 n_2) \, dn_2 \, dn_1,$$

$$(3.10) \qquad J'(w,t) = \int_Z \int_{N \cap w^{-1}\overline{N}w \backslash N} \int_N \Xi'({}^t n_2 w z t n_1) \theta(n_1 n_2) \, dn_2 \, dn_1 dz.$$

Our main result compares these integrals to the integrals $F(w,a)$ introduced in paragraph (2.3).

(3.4) First, for completeness, we record the value of $J(w_0, I_4)$. The proof is similar to that given in paragraph (2.4).

PROPOSITION 3.1.
$$J(w_0, I_4) = 1.$$

To study the remaining integrals, let $\chi_1, \chi_2$ be quasicharacters of $K^\times$. Let us introduce the following formal Mellin transforms of the functions $J'(w,t)$:

(3.11)
$$M_J(\chi_1; w_1) = \int_{K^\times} J'\left(w_1, \mathrm{diag}(a, 1, -1, -a^{-1})\right) \chi_1(a) \, d^\times a,$$

(3.12)
$$M_J(\chi_1; w_2) = \int_{K^\times} J'\left(w_2, \mathrm{diag}(a, a, 1, 1)\right) \chi_1(a) \, d^\times a,$$

(3.13)
$$M_J(\chi_1, \chi_2; 1) = \int_{(K^\times)^2} J'\left(1, \mathrm{diag}(a_1 a_2, a_2, 1, a_1^{-1})\right) \chi_1(a_1) \chi_2(a_2) \, d^\times a_1 \, d^\times a_2.$$

The quantity $M_J(\chi_1, \chi_2; 1)$ is evaluated in Chapter 4 and the other quantities in Chapter 5.

## 4. The Main Theorem

(4.1) To state the main theorem, we let $\sigma a \leftrightarrow wt$ be the bijection between relevant double coset representatives on the two groups described as follows. The double cosets $wt$ for $GSp(4)$ are taken modulo the center. Then the bijection is given by

$$\begin{aligned}(12) \, \mathrm{diag}(a_1 a_2, a_1, 1, 1) &\leftrightarrow \mathrm{diag}(-a_1 a_2, -a_2, 1, a_1^{-1}) \\ \mathrm{diag}(-a_1, a_1, 1, 1) &\leftrightarrow w_1 \, \mathrm{diag}(a_1, 1, -1, -a_1^{-1}) \\ (132) \, \mathrm{diag}(a_2, 1, 1, 1) &\leftrightarrow w_2 \, \mathrm{diag}(a_2, a_2, 1, 1) \\ (1243) &\leftrightarrow w_0.\end{aligned}$$

The main theorem states that the orbital integrals introduced above for the two different groups are in fact equal.

MAIN THEOREM. *For the bijection of relevant cosets $\sigma a \leftrightarrow wt$, the relative orbital integral on $GL(4)$ with respect to the Shalika subgroup is equal to the Kloosterman orbital integral on $GSp(4)$:*

$$F(\sigma, a) = J'(w, t).$$

To establish the Main Theorem, we shall compare these two functions of the form $f(a_1, a_2)$ (resp. $f(a_1)$, $f(a_2)$), by comparing their formal Mellin transforms in the sense described above. As we shall show in this work, these formal Mellin transforms are in fact equal. For the main relevant orbits, these Mellin transforms are computed in Chapters II ($M_F$) and IV ($M_J$), and the comparison is also included in Chapter IV. For the other two nontrivial relevant orbits, the Mellin transforms are computed in Chapters III ($M_F$) and V ($M_J$), and the comparisons are included in Chapter V.

## 5. Kloosterman Integrals and Gauss Sums

(5.1) We end this Chapter with a brief review of Kloosterman integrals and of Gauss sums. For $r, s \in K^\times$, the $(GL(2))$ *Kloosterman integral* $Kl(r, s)$ is defined by the formula

$$Kl(r, s) = \int_{|rx|=1} \psi\left(x - \frac{s}{rx}\right) dx.$$

Here $dx$ denotes the Haar measure on $K$, normalized so that $\mathcal{O}$ has measure 1. The following proposition allows one to evaluate this integral explicitly when $|rs| \neq 1$.

PROPOSITION 5.1. *The Kloosterman integral satisfies the following properties:*
(1) $Kl(r, s) = |rs|^{-1} Kl(s^{-1}, r^{-1})$.
(2) *Suppose $|rs| < 1$.*
    (i) *If $|r| \geq 1$, then $Kl(r, s) = (1 - q^{-1})|r|^{-1}$.*
    (ii) *If $|r| < 1$, $|s| > 1$, then $Kl(r, s) = 0$.*
    (iii) *If $|r| < |\varpi|$, $|s| \leq 1$, then $Kl(r, s) = 0$.*
    (iv) *If $|r| = |\varpi|$, $|s| \leq 1$, then $Kl(r, s) = -1$.*

(5.2) If $\chi$ is a ramified quasicharacter of $K$ and $f$ is a positive integer, we shall say that $\chi$ has *conductor* $f$ if $\chi$ is trivial on the subgroup $1 + \varpi^f \mathcal{O}$ of $\mathcal{O}^\times$, but not on the subgroup $(1 + \varpi^{f-1}\mathcal{O}) \cap \mathcal{O}^\times$. If $\chi$ is a ramified quasicharacter of conductor $f$, let $G(\chi)$ denote the Gauss sum

$$G(\chi) = \int_{|s|=q^f} \chi(s)\psi(s)\, ds.$$

Then one has

(5.1) $$\int_{|s|=q^h} \chi(s)\psi(s)\, ds = \delta_{h,f} G(\chi),$$

where $\delta_{h,f}$ denotes the Kronecker delta, for if $h \neq f$ then the oscillations give zero.